# Additive property of pseudo Drazin inverse of elements in a Banach algebra

Huihui Zhu, Jianlong Chen*

**Abstract:** We study properties of pseudo Drazin inverse in a Banach algebra with unity 1. If $ab = ba$ and $a, b$ are pseudo Drazin invertible, we prove that $a + b$ is pseudo Drazin invertible if and only if $1 + a^{\ddagger}b$ is pseudo Drazin invertible. Moreover, the formula of $(a+b)^{\ddagger}$ is presented . When the commutative condition is weaken to $ab = \lambda ba$ ($\lambda \neq 0$), we also show that $a - b$ is pseudo Drazin invertible if and only if $aa^{\ddagger}(a-b)bb^{\ddagger}$ is pseudo Drazin invertible.

## 1  Introduction

Throughout this paper, $\mathscr{A}$ is a complex Banach algebra with unity 1. The symbols $J(\mathscr{A})$, $\mathscr{A}^{\#}$, $\mathscr{A}^{\text{nil}}$ denote, respectively, the Jacobson radical, the sets of all group invertible, nilpotent elements of $\mathscr{A}$. An element $a \in \mathscr{A}$ is said to have a *Drazin inverse* [10] if there exists $b \in \mathscr{A}$ satisfying

$$ab = ba,\ bab = b,\ a - a^2b \in \mathscr{A}^{\text{nil}}.$$

The element $b$ above is unique and is denoted by $a^D$, and the nilpotency index of $a - a^2b$ is called the Drazin index of $a$, denoted by $\text{ind}(a)$. If $\text{ind}(a) = 1$, then $b$ is the group inverse of $a$ and is denoted by $a^{\#}$. In 2012, Wang and Chen [16] introduced the notion of pseudo Drazin inverse (or p-Drazin inverse for short) in associative rings and Banach algebras. An element $a \in \mathscr{A}$ is called pseudo Drazin invertible if there exists $b \in \mathscr{A}$ such that

$$ab = ba,\ bab = b,\ a^k - a^{k+1}b \in J(\mathscr{A}) \text{ for some integer } k \geqslant 1.$$

*Corresponding author. Department of Mathematics, Southeast University, Nanjing 210096, China. Email: ahzhh08@sina.com (H. Zhu), jlchen@seu.edu.cn (J. Chen)

Any element $b \in \mathscr{A}$ satisfying the conditions above is called a p-Drazin inverse of $a$, denoted by $a^{\ddagger}$. The set of all p-Drazin invertible elements of $\mathscr{A}$ is denoted by $\mathscr{A}^{pD}$. By $a^{\Pi} = 1 - aa^{\ddagger}$ we mean the strongly spectral idempotent of $a$. It is well known that $J(\mathscr{A}) \subset \mathscr{A}^{\text{qnil}}$ and $a \in \mathscr{A}^{\text{qnil}}$ if and only if $\|a^n\|^{\frac{1}{n}} \to 0$ $(n \to \infty)$ in a Banach algebra $\mathscr{A}$.

Wei and Deng [17] presented some additive results on the Drazin inverse for square complex matrices that commute. Zhuang, Chen et al. [18] extended the results in [17] to the ring case and proved that $a + b$ is Drazin invertible if and only if $1 + a^D b$ is Drazin invertible for two commutative Drazin invertible elements $a$ and $b$. Deng [7] explored the Drazin inverse of bounded operators and characterized the relationships between the Drazin inverse of $P - Q$ and $PP^D(P-Q)QQ^D$ with commutativity up to a factor under the condition $PQ = \lambda QP$. More results on (generalized) Drazin inverse can be found in [1-6,8,9,11-15].

Motivated by the paper of Deng [7], Wei, Deng [17] and Zhuang et al. [18], several results in [7,17,18] on Drazin inverse are considered for pseudo Drazin inverse in a Banach algebra. In this paper, we study the product property of pseudo Drazin inverse of two elements and give the equivalent condition of existence on pseudo Drazin inverse of difference of two elements in a Banach algebra under some condition.

## 2 p-Drazin inverse under the condition $ab = ba$

In this section, we give some elementary properties on p-Drazin inverse.

**Lemma 2.1.** [16, Proposition 5.2] *Let $a, b \in \mathscr{A}$. If $ab = ba$ and $a^{\ddagger}$, $b^{\ddagger}$ exist, then $(ab)^{\ddagger} = a^{\ddagger}b^{\ddagger}$.*

In particular, if $a \in \mathscr{A}^{pD}$, then $(a^2)^{\ddagger} = (a^{\ddagger})^2$.

**Lemma 2.2.** *Let $a, b \in \mathscr{A}$. The following statements hold:*
(1) *If $a \in J(\mathscr{A})$, then $ab, ba \in J(\mathscr{A})$,*
(2) *If $a, b \in J(\mathscr{A})$, then $(a+b)^k \in J(\mathscr{A})$ for integer $k \geqslant 1$.*

In [10], some properties of Drazin inverse were presented. One may suspect that if the similar properties can be inherited to p-Drazin inverse. The following results illustrate the possibility.

**Theorem 2.3.** Let $a \in \mathscr{A}^{pD}$. Then

(1) $(a^n)^{\ddagger} = (a^{\ddagger})^n$, $n = 1, 2, \cdots$

(2) $(a^{\ddagger})^{\ddagger} = a^2 a^{\ddagger}$,

(3) $((a^{\ddagger})^{\ddagger})^{\ddagger} = a^{\ddagger}$,

(4) $a^{\ddagger}(a^{\ddagger})^{\ddagger} = aa^{\ddagger}$.

*Proof.* (1) It is obvious when $n = 1$.

Assume the result holds for $n - 1$, i.e., $(a^{n-1})^{\ddagger} = (a^{\ddagger})^{n-1}$.

For $n$, by Lemma 2.1, we have $(a^n)^{\ddagger} = (aa^{n-1})^{\ddagger} = a^{\ddagger}(a^{n-1})^{\ddagger} = a^{\ddagger}(a^{\ddagger})^{n-1} = (a^{\ddagger})^n$.

Hence, $a^n$ is p-Drazin invertible and $(a^n)^{\ddagger} = (a^{\ddagger})^n$.

(2) It is easy to check $a^{\ddagger}a^2 a^{\ddagger} = a^2 a^{\ddagger}a^{\ddagger}$ and $a^2 a^{\ddagger}a^{\ddagger}a^2 a^{\ddagger} = a^2 a^{\ddagger}$.

Since $(a^{\ddagger})^k - (a^{\ddagger})^{k+1}a^2 a^{\ddagger} = (a^{\ddagger})^k - (a^{\ddagger})^{k+1}a = (a^{\ddagger})^k - (a^{\ddagger})^k = 0 \in J(\mathscr{A})$ for some $k \geqslant 1$, it follows that $(a^{\ddagger})^{\ddagger} = a^2 a^{\ddagger}$.

(3) By (2) and Lemma 2.1.

(4) According to (2). □

**Corollary 2.4.** Let $a \in \mathscr{A}^{pD}$. Then $(a^{\ddagger})^{\ddagger} = a$ if and only if $a \in \mathscr{A}^{\#}$.

**Theorem 2.5.** Let $a, b \in \mathscr{A}^{pD}$ with $ab = ba = 0$. Then $(a+b)^{\ddagger} = a^{\ddagger} + b^{\ddagger}$.

*Proof.* Since $ab = ba = 0$, it follows that $ab^{\ddagger} = ba^{\ddagger} = 0$ and $a^{\ddagger}b = b^{\ddagger}a = 0$.

Thus, $(a^{\ddagger} + b^{\ddagger})(a + b) = (a + b)(a^{\ddagger} + b^{\ddagger})$.

(2) Note that $a^{\ddagger}aa^{\ddagger} = a^{\ddagger}$. We have

$$(a^{\ddagger} + b^{\ddagger})(a+b)(a^{\ddagger} + b^{\ddagger}) = (a^{\ddagger}a + b^{\ddagger}b)(a^{\ddagger} + b^{\ddagger})$$
$$= a^{\ddagger} + b^{\ddagger}.$$

(3) According to $a^k - a^{k+1}a^{\ddagger} \in J(\mathscr{A})$ and $b^k - b^{k+1}b^{\ddagger} \in J(\mathscr{A})$ for some $k \geqslant 1$, we obtain

$$(a+b)^k - (a+b)^{k+1}(a^{\ddagger} + b^{\ddagger}) = (a^k + b^k) - (a^{k+1} + b^{k+1})(a^{\ddagger} + b^{\ddagger})$$
$$= a^k + b^k - a^{k+1}a^{\ddagger} - b^{k+1}b^{\ddagger}$$
$$= a^k - a^{k+1}a^{\ddagger} + b^k - b^{k+1}b^{\ddagger}$$
$$\in J(\mathscr{A}).$$

Thus, $(a+b)^{\ddagger} = a^{\ddagger} + b^{\ddagger}$. $\square$

**Corollary 2.6.** *If $a_1, a_2 \cdots a_n$ are p-Drazin invertible of $\mathscr{A}$ with $a_i a_j = 0 (i, j = 1, \cdots, n; i \neq j)$, then $a_1 + a_2 + \cdots + a_n$ is p-Drazin invertible and $(a_1 + a_2 + \cdots + a_n)^{\ddagger} = a_1^{\ddagger} + a_2^{\ddagger} + \cdots + a_n^{\ddagger}$.*

*Proof.* It is right for $n = 2$ by Theorem 2.5.

Assume that the result holds for $n-1$. Then $(a_1 + \cdots + a_{n-1})^{\ddagger} = a_1^{\ddagger} + \cdots + a_{n-1}^{\ddagger}$.

For $n$, since $a_i a_j = 0$, by Theorem 2.5, we have

$$(a_1 + \cdots + a_{n-1} + a_n)^{\ddagger} = (a_1 + \cdots + a_{n-1})^{\ddagger} + a_n^{\ddagger}$$
$$= a_1^{\ddagger} + \cdots + a_n^{\ddagger}.$$

$\square$

In [17], Wei and Deng presented the formula for the Drazin inverse of two square matrices that commute with each other. We consider the result in [17] for p-Drazin inverse in a Banach algebra as follows.

**Theorem 2.7.** *If $a, b \in \mathscr{A}^{pD}$ and $ab = ba$, then $a + b$ is p-Drazin invertible if and only if $1 + a^{\ddagger}b$ is p-Drazin invertible. In this case, we have*

$$(a+b)^{\ddagger} = (1+a^{\ddagger}b)^{\ddagger} a^{\ddagger} + b^{\ddagger} \sum_{i=0}^{\infty} (-b^{\ddagger} a a^{\Pi})^i a^{\Pi},$$

*and $(1 + a^{\ddagger}b)^{\ddagger} = a^{\Pi} + a^2 a^{\ddagger}(a+b)^{\ddagger}$.*

*Proof.* First we show that $1 + b^{\ddagger}aa^{\Pi}$ is invertible. By property of p-Drazin inverse, there exists some $k \geqslant 1$ such that $(aa^{\Pi})^k \in J(\mathscr{A})$. Hence $1 + (b^{\ddagger})^k (aa^{\Pi})^k \in \mathscr{A}^{-1}$. By Lemma 2.1, $a$, $a^{\Pi}$ and $b^{\ddagger}$ commute with each other. We have $(1 + b^{\ddagger}aa^{\Pi})^{-1} = \sum_{i=0}^{\infty}(-b^{\ddagger}aa^{\Pi})^i$.

Suppose that $a + b$ is p-Drazin invertible. We prove that $1 + a^{\ddagger}b$ is p-Drazin invertible. Write $1 + a^{\ddagger}b = a_1 + b_1$ with $a_1 = a^{\Pi}$ and $b_1 = a^{\ddagger}(a+b)$.

Note that $a, b, a^{\ddagger}, b^{\ddagger}$ commute with each other. It follows that $(b_1)^{\ddagger} = (a^{\ddagger}(a+b))^{\ddagger} = a^2 a^{\ddagger}(a+b)^{\ddagger}$ by Lemma 2.1 and Theorem 2.3(2).

Since $a_1$ is idempotent, $(a_1)^{\ddagger} = a^{\Pi}$. Observing that $a_1 b_1 = b_1 a_1 = 0$, it follows that $(1 + a^{\ddagger}b)^{\ddagger} = a^{\Pi} + a^2 a^{\ddagger}(a+b)^{\ddagger}$ from Theorem 2.5.

Conversely, let $x = (1 + a^{\ddagger}b)^{\ddagger}a^{\ddagger} + b^{\ddagger}\sum_{i=0}^{\infty}(-b^{\ddagger}aa^{\Pi})^{i}a^{\Pi}$. We prove $x$ is the Drazin inverse of $a + b$. It is similar to the proof of [18, Theorem 3]. Zhuang, Chen et al. essentially proved that

$$a + b - (a+b)^2 x = a^{\Pi}bb^{\Pi} + aa^{\Pi}(\xi^{\ddagger} - bb^{\ddagger}) + \xi\xi^{\pi}a.$$

We have $(aa^{\Pi})^{k_1} \in J(\mathscr{A})$, $(bb^{\Pi})^{k_2} \in J(\mathscr{A})$ and $(\xi\xi^{\Pi})^{k_3} \in J(\mathscr{A})$ for some $k_1$, $k_2$ and $k_3$. Take suitable $k \geqslant \max\{k_1 + k_2 + k_3\}$, it follows that $(a + b)^k - (a + b)^{k+1}x = (a + b - (a+b)^2 x)^k \in J(\mathscr{A})$ by Lemma 2.2(2).

The proof is completed. □

**Corollary 2.8.** *Let $a, b \in \mathscr{A}^{pD}$ with $ab = ba$. If $1 + a^{\ddagger}b$ is p-Drazin invertible, then*

(1) *If $a$ is nilpotent, then $(a + b)^{\ddagger} = b^{\ddagger}$,*

(2) *If $a$ is invertible, then $(a + b)^{\ddagger} = (1 + a^{-1}b)^{\ddagger}a^{-1} + b^{\ddagger}a^{-1}$,*

(3) *If $a$ is group invertible, then $(a + b)^{\ddagger} = (1 + a^{\#}b)^{\ddagger}a^{\#} + a^{\Pi}b^{\ddagger}$.*

## 3  p-Drazin inverse under the condition $ab = \lambda ba$

In this section, we give some results on p-Drazin inverse under the condition that $ab = \lambda ba$.

**Lemma 3.1.** *Let $a, b \in \mathscr{A}$ with $ab = \lambda ba$ $(\lambda \neq 0)$. Then*

(1) $ab^n = \lambda^n b^n a$, $a^n b = \lambda^n b a^n$,

(2) $(ab)^n = \lambda^{\frac{-n(n-1)}{2}} a^n b^n = \lambda^{\frac{n(n+1)}{2}} b^n a^n$,

(3) $(ba)^n = \lambda^{\frac{n(n-1)}{2}} b^n a^n = \lambda^{\frac{-n(n+1)}{2}} a^n b^n$.

*Proof.* By induction, it is easy to obtain the results. □

It is well known that $c^k \in J(\mathscr{A})$ implies that $c \in \mathscr{A}^{\text{qnil}}$ for $k \geqslant 1$. For any $x \in \text{comm}(c)$, $1 - (cx)^k = (1 - cx)(1 + cx + \cdots + (cx)^{k-1}) \in \mathscr{A}^{-1}$. Hence, $1 - cx \in \mathscr{A}^{-1}$ implies that $c \in \mathscr{A}^{\text{qnil}}$. By property of p-Drazin inverse, we have $a(1 - aa^{\ddagger}) \in \mathscr{A}^{\text{qnil}}$

**Lemma 3.2.** *Let $a, b \in \mathscr{A}^{pD}$ with $ab = \lambda ba$ $(\lambda \neq 0)$. Then*

(1) $aa^{\ddagger}b = baa^{\ddagger}$,

(2) $bb^{\ddagger}a = abb^{\ddagger}$.

*Proof.* The Proof essentially in [4]. Since $a(1 - aa^{\ddagger}) \in \mathscr{A}^{\text{qnil}}$, it follows that $\|(a(1 - aa^{\ddagger}))^n\|^{\frac{1}{n}} \to 0$ $(n \to \infty)$. Suppose $p = aa^{\ddagger}$. We have

$$\begin{aligned}
\|pb - pbp\|^{\frac{1}{n}} &= \|a^{\ddagger}ab(1 - aa^{\ddagger})\|^{\frac{1}{n}} \\
&= \|(a^{\ddagger})^n a^n b(1 - aa^{\ddagger})\|^{\frac{1}{n}} \\
&= \|(a^{\ddagger})^n \lambda^n ba^n (1 - aa^{\ddagger})\|^{\frac{1}{n}} \\
&= \|\lambda^n (a^{\ddagger})^n b(a(1 - aa^{\ddagger}))^n\|^{\frac{1}{n}} \\
&\leqslant |\lambda| \, \|a^{\ddagger}\| \|b\|^{\frac{1}{n}} \|(a(1 - aa^{\ddagger}))^n\|^{\frac{1}{n}}.
\end{aligned}$$

Hence, $\|pb - pbp\|^{\frac{1}{n}} \to 0$ $(n \to \infty)$, it follows that $pb = pbp$.

Similarly, $bp = pbp$.

Thus, $aa^{\ddagger}b = baa^{\ddagger}$.

(2) The proof is similar (1). $\square$

Wang and Chen [16] proved that $(ab)^{\ddagger} = b^{\ddagger}a^{\ddagger}$ under the condition $ab = ba$. We obtain some generalized results under weak commutative condition $ab = \lambda ba$.

**Theorem 3.3.** *Let $a, b \in \mathscr{A}^{pD}$ with $ab = \lambda ba$ $(\lambda \neq 0)$. Then*
(1) $a^{\ddagger}b = \lambda^{-1}ba^{\ddagger}$,
(2) $ab^{\ddagger} = \lambda^{-1}b^{\ddagger}a$,
(3) $(ab)^{\ddagger} = b^{\ddagger}a^{\ddagger} = \lambda^{-1}a^{\ddagger}b^{\ddagger}$.

*Proof.* (1) Since $aa^{\ddagger}b = baa^{\ddagger}$, we have

$$\begin{aligned}
a^{\ddagger}b &= a^{\ddagger}aa^{\ddagger}b = a^{\ddagger}baa^{\ddagger} = a^{\ddagger}\lambda^{-1}aba^{\ddagger} \\
&= \lambda^{-1}a^{\ddagger}aba^{\ddagger} = \lambda^{-1}baa^{\ddagger}a^{\ddagger} \\
&= \lambda^{-1}ba^{\ddagger}.
\end{aligned}$$

(2) The proof is similar to (1).

(3) We first prove that $b^{\ddagger}a^{\ddagger} = \lambda^{-1}a^{\ddagger}b^{\ddagger}$. By Lemma 3.2, it follows that

$$\begin{aligned}
b^{\ddagger}a^{\ddagger} &= b^{\ddagger}(aa^{\ddagger})a^{\ddagger} = (aa^{\ddagger})b^{\ddagger}a^{\ddagger} \\
&= a^{\dagger}(ab^{\ddagger})a^{\ddagger} = a^{\ddagger}(\lambda^{-1}b^{\ddagger}a)a^{\ddagger} \\
&= \lambda^{-1}a^{\ddagger}b^{\ddagger}(aa^{\ddagger}) = \lambda^{-1}a^{\ddagger}(aa^{\ddagger})b^{\ddagger}
\end{aligned}$$

$$= \lambda^{-1}a^{\ddagger}b^{\ddagger}.$$

Let $x = b^{\ddagger}a^{\ddagger}$. We prove $x$ is the p-Drazin inverse of $ab$.

(1) By Lemma 3.2, we obtain

$$\begin{aligned}(ab)x &= abb^{\ddagger}a^{\ddagger} = aa^{\ddagger}b^{\ddagger}b \\ &= b^{\ddagger}aa^{\ddagger}b = b^{\ddagger}a^{\ddagger}ab \\ &= x(ab).\end{aligned}$$

(2) $x(ab)x = b^{\ddagger}(a^{\ddagger}a)bb^{\ddagger}a^{\ddagger} = b^{\ddagger}bb^{\ddagger}(a^{\ddagger}a)a^{\ddagger} = b^{\ddagger}a^{\ddagger} = x.$

(3) We present an useful equality, i.e.,

$$\begin{aligned}b^{k+1}b^{\ddagger}a^{k+1}a^{\ddagger} &= b^{k}bb^{\ddagger}a^{k}aa^{\ddagger} = b^{k}a^{k}aa^{\ddagger}bb^{\ddagger} \\ &= b^{k}a^{k+1}a^{\ddagger}bb^{\ddagger} = b^{k}a^{k+1}bb^{\ddagger}a^{\ddagger} \\ &= b^{k}\lambda^{k+1}ba^{k+1}b^{\ddagger}a^{\ddagger} \\ &= \lambda^{k+1}b^{k+1}a^{k+1}b^{\ddagger}a^{\ddagger}.\end{aligned}$$

Hence, we have

$$\begin{aligned}(ab)^{k} - (ab)^{k+1}b^{\ddagger}a^{\ddagger} &= \lambda^{\frac{k(k+1)}{2}}b^{k}a^{k} - \lambda^{\frac{(k+1)(k+2)}{2}}b^{k+1}a^{k+1}b^{\ddagger}a^{\ddagger} \\ &= \lambda^{\frac{k(k+1)}{2}}(b^{k}a^{k} - \lambda^{k+1}b^{k+1}a^{k+1}b^{\ddagger}a^{\ddagger}) \\ &= \lambda^{\frac{k(k+1)}{2}}(b^{k}a^{k} - b^{k+1}b^{\ddagger}a^{k+1}a^{\ddagger}) \\ &= \lambda^{\frac{k(k+1)}{2}}[-(b^{k} - b^{k+1}b^{\ddagger})(a^{k} - a^{k+1}a^{\ddagger}) + (b^{k} - b^{k+1}b^{\ddagger})a^{k} \\ &\quad + b^{k}(a^{k} - a^{k+1}a^{\ddagger})] \\ &\in J(\mathscr{A})\end{aligned}$$

for some $k \geqslant 1$.

Therefore, $(ab)^{\ddagger} = b^{\ddagger}a^{\ddagger} = \lambda^{-1}a^{\ddagger}b^{\ddagger}$. $\square$

**Corollary 3.4.** *Let $a, b \in \mathscr{A}^{pD}$ with $ab = \lambda ba$ ($\lambda \neq 0$). Then*
(1) $(a^{\ddagger}b)^{n} = \lambda^{\frac{n(n-1)}{2}}(a^{\ddagger})^{n}b^{n} = \lambda^{\frac{-n(n+1)}{2}}b^{n}(a^{\ddagger})^{n}$,
(2) $(ab^{\ddagger})^{n} = \lambda^{\frac{n(n-1)}{2}}a^{n}(b^{\ddagger})^{n} = \lambda^{\frac{-n(n+1)}{2}}(b^{\ddagger})^{n}a^{n}$.

*Proof.* By induction and Theorem 3.3. $\square$

**Theorem 3.5.** Let $a, b \in \mathscr{A}^{pD}$ with $ab = \lambda ba$ ($\lambda \neq 0$). Then $a - b$ is p-Drazin invertible if and only if $w = aa^{\ddagger}(a - b)bb^{\ddagger}$ is p-Drazin invertible. In this case, we have

$$(a - b)^{\ddagger} = w^{\ddagger} + a^{\ddagger} \sum_{i=0}^{\infty} (ba^{\ddagger})^i b^{\Pi} - a^{\Pi} \sum_{i=0}^{\infty} (b^{\ddagger}a)^i b^{\ddagger},$$

and $w^{\ddagger} = a^{\ddagger}(a - b)^{\ddagger}bb^{\ddagger}$.

*Proof.* Assume that $a - b$ is p-Drazin invertible. Since $aa^{\ddagger}$ is idempotent, $aa^{\ddagger} \in \mathscr{A}^{pD}$. By Lemma 3.2, we have $aa^{\ddagger}(a - b) = (a - b)aa^{\ddagger}$. Hence, $aa^{\ddagger}(a - b) \in \mathscr{A}^{pD}$ according to Lemma 2.1.

From 3.2, it follows that $aa^{\ddagger}(a - b)bb^{\ddagger} = bb^{\ddagger}aa^{\ddagger}(a - b)$. Thus, by Lemma 2.1, $aa^{\ddagger}(a - b)bb^{\ddagger} \in \mathscr{A}^{pD}$ and $w^{\ddagger} = a^{\ddagger}(a - b)^{\ddagger}bb^{\ddagger}$.

Conversely, let

$$x = w^{\ddagger} + a^{\ddagger} \sum_{i=0}^{\infty} (ba^{\ddagger})^i b^{\Pi} - a^{\Pi} \sum_{i=0}^{\infty} (b^{\ddagger}a)^i b^{\ddagger}.$$

We prove that $x$ is the p-Drazin inverse of $a - b$.

(1) By Lemma 3.2, we have

$$\begin{aligned}
(a - b)w &= (a - b)aa^{\ddagger}(a - b)bb^{\ddagger} \\
&= aa^{\ddagger}(a - b)(a - b)bb^{\ddagger} \\
&= aa^{\ddagger}(a - b)bb^{\ddagger}(a - b) \\
&= w(a - b).
\end{aligned}$$

Hence, $(a - b)w^{\ddagger} = w^{\ddagger}(a - b)$. Moreover,

$$\begin{aligned}
(a - b)a^{\ddagger} \sum_{i=0}^{\infty} (ba^{\ddagger})^i b^{\Pi} &= (aa^{\ddagger} - ba^{\ddagger}) \sum_{i=0}^{\infty} (ba^{\ddagger})^i b^{\Pi} = aa^{\ddagger} \sum_{i=0}^{\infty} (ba^{\ddagger})^i b^{\Pi} - ba^{\ddagger} \sum_{i=0}^{\infty} (ba^{\ddagger})^i b^{\Pi} \\
&= aa^{\ddagger}(1 + ba^{\ddagger} + (ba^{\ddagger})^2 + \cdots)b^{\Pi} - ba^{\ddagger}(1 + ba^{\ddagger} + (ba^{\ddagger})^2 + \cdots)b^{\Pi} \\
&= (aa^{\ddagger} + ba^{\ddagger} + (ba^{\ddagger})^2 + \cdots)b^{\Pi} - (ba^{\ddagger} + (ba^{\ddagger})^2 + \cdots)b^{\Pi} \\
&= aa^{\ddagger}b^{\pi}
\end{aligned}$$

Similarly, $(a - b)(a^{\Pi} \sum_{i=0}^{\infty} (b^{\ddagger}a)^i b^{\ddagger}) = -bb^{\ddagger}a^{\Pi}$.

We have

$$(a - b)x = (a - b)(w^{\ddagger} + a^{\ddagger} \sum_{i=0}^{\infty} (ba^{\ddagger})^i b^{\Pi} - a^{\Pi} \sum_{i=0}^{\infty} (b^{\ddagger}a)^i b^{\ddagger})$$

$$= (a-b)w^{\ddagger} + aa^{\ddagger}b^{\Pi} + bb^{\ddagger}a^{\Pi}.$$

Note that Lemma 3.2. We get

$$a^{\ddagger}\sum_{i=0}^{\infty}(ba^{\ddagger})^{i}b^{\Pi}(a-b) = a^{\ddagger}\sum_{i=0}^{\infty}(ba^{\ddagger})^{i}(a-b)b^{\Pi} = a^{\ddagger}\sum_{i=0}^{\infty}(ba^{\ddagger})^{i}ab^{\Pi} - a^{\ddagger}\sum_{i=0}^{\infty}(ba^{\ddagger})^{i}bb^{\Pi}$$
$$= a^{\ddagger}(1 + ba^{\ddagger} + (ba^{\ddagger})^{2} + \cdots)ab^{\Pi} - a^{\ddagger}(1 + ba^{\ddagger} + (ba^{\ddagger})^{2} + \cdots)bb^{\Pi}$$
$$= a^{\ddagger}(ab^{\Pi} + ba^{\ddagger}ab^{\Pi} + (ba^{\ddagger})^{2}ab^{\Pi} + \cdots) - a^{\ddagger}(bb^{\Pi} + ba^{\ddagger}bb^{\Pi} + (ba^{\ddagger})^{2}bb^{\Pi} + \cdots)$$
$$= a^{\ddagger}(ab^{\Pi} + a^{\ddagger}abb^{\Pi} + aa^{\ddagger}ba^{\ddagger}bb^{\Pi} + aa^{\ddagger}(ba^{\ddagger})^{2}bb^{\Pi} + \cdots)$$
$$\quad - a^{\ddagger}(bb^{\Pi} + ba^{\ddagger}bb^{\Pi} + (ba^{\ddagger})^{2}bb^{\Pi} + \cdots)$$
$$= a^{\ddagger}(ab^{\Pi} + bb^{\Pi} + ba^{\ddagger}bb^{\Pi} + (ba^{\ddagger})^{2}bb^{\Pi} + \cdots)$$
$$\quad - a^{\ddagger}(bb^{\Pi} + ba^{\ddagger}bb^{\Pi} + (ba^{\ddagger})^{2}bb^{\Pi} + \cdots)$$
$$= aa^{\ddagger}b^{\Pi}.$$

Similarly, $a^{\Pi}\sum_{i=0}^{\infty}(b^{\ddagger}a)^{i}b^{\ddagger}(a-b) = -bb^{\ddagger}a^{\Pi}$.

Therefore, we have

$$x(a-b) = (w^{\ddagger} + a^{\ddagger}\sum_{i=0}^{\infty}(ba^{\ddagger})^{i}b^{\Pi} - a^{\Pi}\sum_{i=0}^{\infty}(b^{\ddagger}a)^{i}b^{\ddagger})(a-b)$$
$$= (a-b)w^{\ddagger} + aa^{\ddagger}b^{\Pi} + bb^{\ddagger}a^{\Pi}$$
$$= (a-b)x.$$

(2) From Lemma 2.1, we obtain that $w^{\ddagger} = aa^{\ddagger}(a-b)^{\ddagger}bb^{\ddagger}$ and

$$w^{\ddagger}(a-b)w^{\ddagger} = aa^{\ddagger}(a-b)^{\ddagger}bb^{\ddagger}(a-b)aa^{\ddagger}(a-b)^{\ddagger}bb^{\ddagger}$$
$$= aa^{\ddagger}bb^{\ddagger}(a-b)^{\ddagger}$$
$$= w^{\ddagger}.$$

By property of p-Drazin inverse, we have $w^{\ddagger}a^{\Pi} = a^{\Pi}w^{\ddagger} = 0$ and $w^{\ddagger}b^{\Pi} = b^{\Pi}w^{\ddagger} = 0$.
Hence, $w^{\ddagger}(aa^{\ddagger}b^{\Pi} + bb^{\ddagger}a^{\Pi}) = w^{\ddagger}b^{\Pi}aa^{\ddagger} + w^{\ddagger}a^{\Pi}bb^{\ddagger} = 0$.

According to Lemma 3.2 and $(a-b)w^{\ddagger} = w^{\ddagger}(a-b)$, we have

$$[a^{\ddagger}\sum_{i=0}^{\infty}(ba^{\ddagger})^{i}b^{\Pi} - a^{\Pi}\sum_{i=0}^{\infty}(b^{\ddagger}a)^{i}b^{\ddagger}](a-b)w^{\ddagger} = [a^{\ddagger}\sum_{i=0}^{\infty}(ba^{\ddagger})^{i}b^{\Pi} - a^{\Pi}\sum_{i=0}^{\infty}(b^{\ddagger}a)^{i}b^{\ddagger}]w^{\ddagger}(a-b)$$
$$= 0$$

and

$$[a^{\ddagger}\sum_{i=0}^{\infty}(ba^{\ddagger})^{i}b^{\Pi} - a^{\Pi}\sum_{i=0}^{\infty}(b^{\ddagger}a)^{i}b^{\ddagger}](aa^{\ddagger}b^{\Pi} + bb^{\ddagger}a^{\Pi}) = a^{\ddagger}\sum_{i=0}^{\infty}(ba^{\ddagger})^{i}b^{\Pi}aa^{\ddagger}b^{\Pi}$$

$$+a^{\ddagger}\sum_{i=0}^{\infty}(ba^{\ddagger})^{i}b^{\Pi}bb^{\ddagger}a^{\Pi} - a^{\Pi}\sum_{i=0}^{\infty}(b^{\ddagger}a)^{i}b^{\ddagger}aa^{\ddagger}b^{\Pi} - a^{\Pi}\sum_{i=0}^{\infty}(b^{\ddagger}a)^{i}b^{\ddagger}bb^{\ddagger}a^{\Pi}$$

$$= a^{\ddagger}\sum_{i=0}^{\infty}(ba^{\ddagger})^{i}b^{\Pi} - a^{\Pi}\sum_{i=0}^{\infty}(b^{\ddagger}a)^{i}b^{\ddagger}.$$

Therefore,

$$x(a-b)x = [w^{\ddagger} + a^{\ddagger}\sum_{i=0}^{\infty}(ba^{\ddagger})^{i}b^{\Pi} - a^{\Pi}\sum_{i=0}^{\infty}(b^{\ddagger}a)^{i}b^{\ddagger}][(a-b)w^{\ddagger} + aa^{\ddagger}b^{\Pi} + bb^{\ddagger}a^{\Pi}]$$

$$= w^{\ddagger}(a-b)w^{\ddagger} + w^{\ddagger}(aa^{\ddagger}b^{\Pi} + bb^{\ddagger}a^{\Pi}) + [a^{\ddagger}\sum_{i=0}^{\infty}(ba^{\ddagger})^{i}b^{\Pi} - a^{\Pi}\sum_{i=0}^{\infty}(b^{\ddagger}a)^{i}b^{\ddagger}](a-b)w^{\ddagger}$$

$$+[a^{\ddagger}\sum_{i=0}^{\infty}(ba^{\ddagger})^{i}b^{\Pi} - a^{\Pi}\sum_{i=0}^{\infty}(b^{\ddagger}a)^{i}b^{\ddagger}](aa^{\ddagger}b^{\Pi} + bb^{\ddagger}a^{\Pi})$$

$$= w^{\ddagger} + [a^{\ddagger}\sum_{i=0}^{\infty}(ba^{\ddagger})^{i}b^{\Pi} - a^{\Pi}\sum_{i=0}^{\infty}(b^{\ddagger}a)^{i}b^{\ddagger}](aa^{\ddagger}b^{\Pi} + bb^{\ddagger}a^{\Pi})$$

$$= w^{\ddagger} + a^{\ddagger}\sum_{i=0}^{\infty}(ba^{\ddagger})^{i}b^{\Pi} - a^{\Pi}\sum_{i=0}^{\infty}(b^{\ddagger}a)^{i}b^{\ddagger}$$

$$= x.$$

(3) Since $w^{\ddagger} = aa^{\ddagger}(a-b)^{\ddagger}bb^{\ddagger}$, we have

$$(a-b)^{2}w^{\ddagger} = (a-b)(a-b)aa^{\ddagger}(a-b)^{\ddagger}bb^{\ddagger}$$

$$= aa^{\ddagger}(a-b)^{\ddagger}bb^{\ddagger}aa^{\ddagger}(a-b)bb^{\ddagger}aa^{\ddagger}(a-b)bb^{\ddagger}$$

$$= w^{2}w^{\ddagger} = w(1-w^{\Pi})$$

$$= w - ww^{\Pi},$$

and

$$(a-b)(aa^{\ddagger}b^{\Pi} + bb^{\ddagger}a^{\Pi}) = (a-b)(1-a^{\Pi})b^{\Pi} + (a-b)(1-b^{\Pi})a^{\Pi}$$

$$= ab^{\Pi} - ba^{\Pi} + aa^{\Pi} - bb^{\Pi} - 2aa^{\Pi}b^{\Pi} + 2a^{\Pi}bb^{\Pi}.$$

Observing that $w = (1-a^{\Pi})(a-b)(1-bb^{\Pi}) = (a-b)-(a-b)bb^{\Pi}-a^{\Pi}(a-b)+a^{\Pi}(a-b)b^{\Pi}$, it follows that

$$(a-b)^{k} - (a-b)^{k+1}x = ((a-b) - (a-b)^{2}x)^{k}$$

$$
\begin{aligned}
&= ((a-b) - (a-b)[(a-b)w^{\ddagger} + aa^{\ddagger}b^{\Pi} + bb^{\ddagger}a^{\Pi}])^k \\
&= ((a-b) - (a-b)^2 w^{\ddagger} - (a-b)(aa^{\ddagger}b^{\Pi} + bb^{\ddagger}a^{\Pi}))^k \\
&= ((a-b) - w + ww^{\ddagger} - (a-b)(aa^{\ddagger}b^{\Pi} + bb^{\ddagger}a^{\Pi}))^k \\
&= ((a-b) - (a-b) + (a-b)bb^{\Pi} + a^{\Pi}(a-b) \\
&\quad + a^{\Pi}(a-b)bb^{\Pi} - ab^{\Pi} + ba^{\Pi} - aa^{\Pi} + bb^{\Pi} + 2aa^{\Pi}b^{\Pi} - 2a^{\Pi}bb^{\Pi})^k \\
&= (bb^{\Pi}a^{\Pi} - aa^{\Pi}b^{\Pi} - ww^{\Pi})^k.
\end{aligned}
$$

By the property of pseudo Drazin inverse, there exist some $k_1$, $k_2$ and $k_3$ such that $(aa^{\Pi})^{k_1} \in J(\mathscr{A})$, $(bb^{\Pi})^{k_2} \in J(\mathscr{A})$ and $(ww^{\Pi})^{k_3} \in J(\mathscr{A})$. Take suitable large $k \geqslant \max\{k_1 + k_2 + k_3\}$, we obtain that $(a-b)^k - (a-b)^{k+1}x \in J(\mathscr{A})$ by Lemma 2.2(2). $\square$

It is well-known that $a \in \mathcal{A}^{\text{nil}}$ implies that $a^k = 0$ for some positive integer $k$. Hence, if $a - a^2 b \in \mathcal{A}^{\text{nil}}$, then $0 = (a - a^2 b)^k \in J(\mathcal{A})$. Let $\mathscr{B}(X)$ be a Banach algebra of all bounded linear operators on $X$. Then we have the following result.

**Corollary 3.6.** [7, Theorem 2.4] *Let $a, b \in \mathscr{B}(X)$ be Drazin invertible with $s = \text{ind}(a)$ and $t = \text{ind}(b)$. If $ab = \lambda ba$ and $\lambda \neq 0$, then $a - b$ is Drazin invertible if and only if $w = aa^D(a-b)bb^D$ is Drazin invertible. In this case, we have*

$$
(a-b)^D = w^D + (1 - bb^D)\left[\sum_{i=0}^{t-1} \lambda^{\frac{i(i+1)}{2}} (a^D)^{i+1} b^i\right] - \left[\sum_{i=0}^{s-1} \lambda^{\frac{i(i+1)}{2}} a^i (b^D)^{i+1}\right](1 - aa^D)
$$

*Proof.* By Corollary 3.4 and Theorem 3.5. $\square$

### ACKNOWLEDGMENTS

This research is supported by the National Natural Science Foundation of China (10971024), the Specialized Research Fund for the Doctoral Program of Higher Education (20120092110020), the Natural Science Foundation of Jiangsu Province (BK2010393) and the Foundation of Graduate Innovation Program of Jiangsu Province(CXLX13-072).

# References

[1] N. Castro-González, Additive perturbation results for the Drazin inverse, Linear Algebra Appl. 397 (2005), 279-297.

[2] N. Castro-González, J.J. Koliha, New additive results for the g-Drazin inverse, Proc. Roy. Soc. Edinburgh Sect. A 134 (2004), 1085-1097.

[3] N. Castro-González, M.F. Martinez-Serrano, Expressions for the g-Drazin inverse of additive perturbed elements in a Banach algebra, Linear Algebra Appl. 432 (2010), 1885-1895.

[4] D.S. Cvetković-Ilić, The generalized Drazin inverse with commutativity up to a factor in a Banach algebra, Linear Algebra Appl. 431 (2009), 783-791.

[5] D.S. Cvetković-Ilić, D.S. Djordjević and Y.M. Wei, Additive results for the generalized Drazin inverse in a Banach algebra, Linear Algebra Appl. 418 (2006), 53-61.

[6] D.S. Cvetković-Ilić, X.J. Liu and Y.M. Wei, Some additive results for the generalized Drazin inverse in Banach algebra, Elect. J. Linear Algebra. 22 (2011), 1049-1058.

[7] C.Y. Deng, The Drazin inverse of bounded operators with commutativity up to a factor, Appl. Math. Comput. 206 (2008), 695-703.

[8] C.Y. Deng, Y.M. Wei, New additive results for the generalized Drazin inverse, J. Math. Anal. Appl. 370 (2010), 313-321.

[9] D.S. Djordjević, Y.M. Wei, Additive results for the generalized Drazin inverse, J. Austral. Math. Soc. 73 (2002), 115-125.

[10] M.P. Drazin, Pseudo-inverses in associative rings and semigroups, Amer. Math. Monthly 65 (1958) 506-514.

[11] R.E. Hartwig, G.R. Wang and Y.M. Wei, Some additive results on Drazin inverse, Linear Algebra Appl. 322 (2001), 207-217.

[12] J.J. Koliha, A generalized Drazin inverse, Glasgow Math. J. 38 (1996), 367-381.
12